\documentclass[invmat]{svjour}
\usepackage{latexsym,amsfonts,amssymb,amsbsy,amsmath}
\usepackage{graphicx}
\usepackage{enumerate}

\def\lover#1#2{\lower #1 ex\hbox{$#2$}}

\def\DJ{{D\kern-.8em\raise.27ex\hbox{-$\!$-}\kern.3em}}

\def\Q{{\mathbb Q}}

\def\nastavak{eps}

\begin{document}

\dedication{Dedicated to the memories of C.\ F.\ Gauss (1777--1855).}

\title{Geometry of pentagons:\\
from Gauss to Robbins}
\titlerunning{Geometry of pentagons}
\authorrunning{Dragutin Svrtan et al.}

\author{Dragutin Svrtan\inst{1}
\and Darko Veljan\inst{2}
\and Vladimir Volenec\inst{3}
}

\institute{Department of Mathematics, University of Zagreb,
Bijeni\v{c}ka cesta 30, 10000 Zagreb, Croatia, \email{dsvrtan@math.hr}
\and
Department of Mathematics, University of Zagreb,
Bijeni\v{c}ka cesta 30, 10000 Zagreb, Croatia, \email{dveljan@math.hr}
\and
Department of Mathematics, University of Zagreb,
Bijeni\v{c}ka cesta 30, 10000 Zagreb, Croatia, \email{volenec@math.hr}
}

\offprints{Dragutin Svrtan}
\mail{dsvrtan@math.hr}

\date{Received: date / Revised version: date}

\maketitle

\begin{abstract}
An almost forgotten gem of Gauss tells us how to compute the area of a
pentagon by just going around it and measuring areas of each vertex
triangles (i.e.\ triangles whose vertices are three consecutive vertices
of the pentagon). We give several proofs and extensions of this beautiful
formula to hexagon etc.\ and consider special cases of affine--regular
polygons.

The Gauss pentagon formula is, in fact, equivalent to the Monge formula which is
equivalent to the Ptolemy formula.

On the other hand, we give a new proof of the Robbins formula for the area of a
cyclic pentagon in terms of the side lengths, and this is a consequence of the Ptolemy
formula. The main tool is simple: just eliminate from algebraic equations, via resultants.
By combining Gauss and Robbins formulas we get an explicit rational expression
for the area of any cyclic pentagon.
So, after centuries of geometry of triangles and quadrilaterals, we arrive to the nontrivial
geometry of pentagons.

\keywords{Gauss pentagon formula, Robbins pentagon formula}

\end{abstract}

{\bf AMS 2000 Mathematics Subject Classification:} 51M04, 51M25,
52A10, 52A38, 52B55

\section{INTRODUCTION}
\label{intro}

As it is well known, the great Gauss initiated not only big mathematical
theories, but also was an author of many small mathematical "gems" in
number theory, arithmetic and algebra, real and complex analysis,
probability theory, classical mechanics and particularly in geometry. One
of his beloved topics in elementary geometry were ordinary plane polygons.
Not only their constructions by ruler and compass, but also computing
their areas. In this paper we discuss one of his less known and almost
forgotten "gems": the pentagon formula. We have not been aware of this
formula until recently, when one of us (D.S.) bought Wu's book \cite{1},
while being at the World Congress of Mathematicians in Beijing in August
2002. Wu in his book on p.\ 337 gives a mechanical (coordinate) proof of
the Gauss pentagon formula in the style of symbolic computations. As we
have been able to trace down, the only other place where this formula is
stated and proved is Gauss' original paper \cite{2} and reproduced in his
Collected works \cite{3}. Occasionally though, some authors discuss
closely related problems, e.g.\ \cite{4}, or also closely related Monge
formula, e.g.\ \cite{5}.

It this paper we offer several proofs of this beautiful formula for the
area of a pentagon and show that, in fact, it is equivalent to the Monge
formula, which, in turn, is equivalent to the well known Ptolemy formula
$ef=ac+bd$, where $a, b, c, d$ are side lengths and $e, f$ diagonal
lengths of a cyclic (inscribed) quadrilateral.

In words, the Gauss pentagon formula says that to compute the area of a
(convex) pentagon we don't have to dissect, triangulate or integrate it,
but only go around the pentagon and measure areas of its vertex triangles
(i.e.\ triangles with three consecutive vertices). We shall also discuss
extensions to hexagons etc.\ and consider affine--regular polygons.

Our main result is the explicit Robbins formula \cite{6} and \cite{77} for
the area of a cyclic pentagon (i.e.\ a pentagon inscribed in a circle) in
terms of its side lengths. It is a consequence of the Ptolemy formula, and
that's the reason why we included it here. In this way we continue the
sequence: the Heron formula (1.\ century), the Brahmagupta formula (7.\
century), the Robbins formula (1995),... The method we use is simple:
eliminate via resultants.

\section{THE GAUSS PENTAGON FORMULA}
\label{sec:1}

\indent
Consider any convex pentagon with vertices (in cyclic order) $0$,
$1$, $2$, $3$, $4$. Denote by $A$ the area of the whole pentagon, and by
$(ijk)$ the area of the triangle $\triangle ijk$, $0\leq i,j,k\leq 4$. In
particular, the area of the vertex triangle $\triangle (i-1,i,i+1)
(\mbox{mod } 5)$, we denote simply by $(i)$.

Next, consider the cyclic symmetric functions $c_1$, $c_2$ of
the first and second degree in variables $(0)$, $(1)$, $(2)$,
$(3)$, $(4)$:
\[
\begin{array}{l}
  c_1 = (0) + (1) + (2) + (3) + (4), \\[2mm]
  c_2 = (0)(1) + (1)(2) + (2)(3) + (3)(4) + (4)(0).
\end{array}
\]
\begin{figure}
\centerline{\includegraphics[width=145pt,height=139pt]{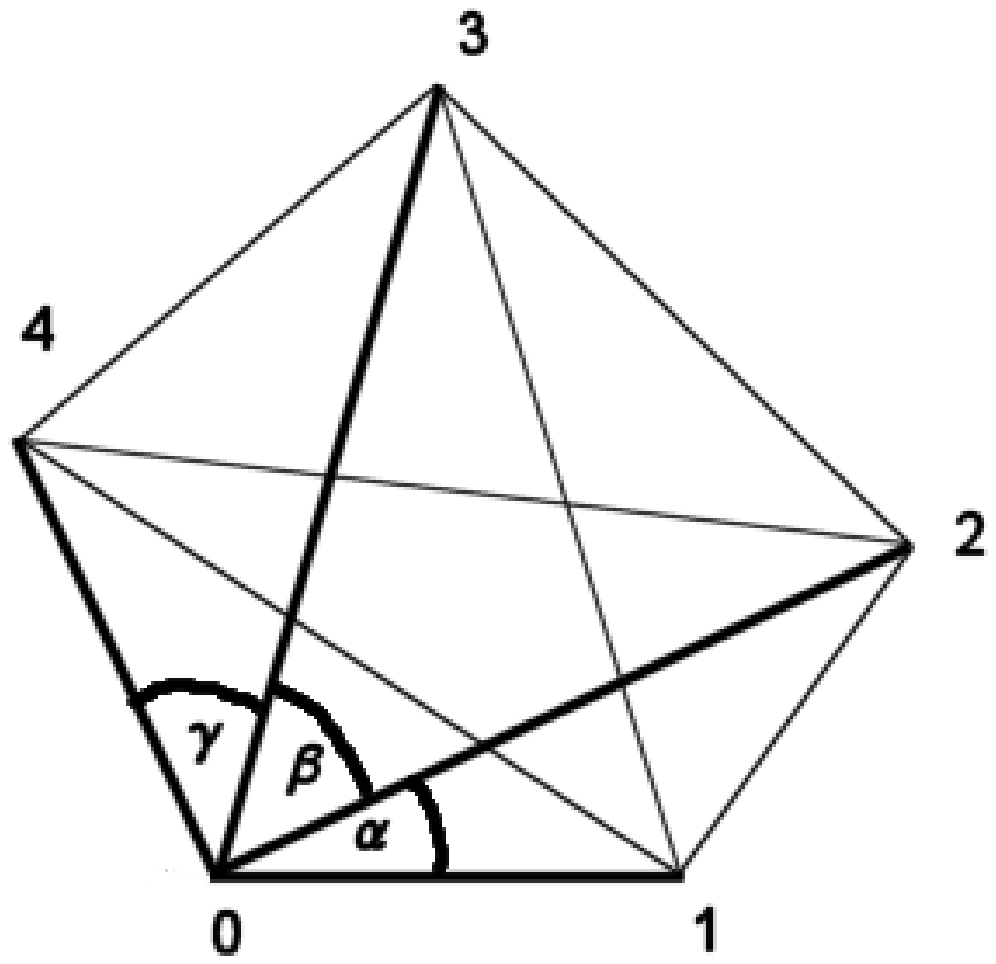}}
\caption{\label{Fig:1}}
\end{figure}
\begin{theorem}
\label{T:Gpf} {\bf (the Gauss pentagon formula)} With above notations we have
\begin{equation}
\label{F:G}\tag{$\ast$}
A^2 -c_1A+c_2=0.
\end{equation}
\end{theorem}
\proof The crucial part in proving the Gauss formula is the following Monge formula (sometimes called Pl\"
ucker's formula) involving areas of all six triangles having $0$ as a vertex:
\begin{equation}
\label{F:M}\tag{$+$}
(012)(034)+(014)(023)=(013)(024).
\end{equation}
Now, the Gauss formula (\ref{F:G}) follows simply from (\ref{F:M}) by just noting from Figure \ref{Fig:1} that
$(023) = A - (1) - (4)$, $(013) = A - (2)-(4)$, $(024)=A-(1)-(3)$ and writing $(014)=(0)$, $(012)=(1)$, and
$(034)=(4)$.

It remains only to prove the Monge formula (\ref{F:M}). By writing each
triangle area as one half of the product of two side lengths and the
sine of the angle between them, we see again from the Figure \ref{Fig:1}
that, after dividing by the product of the side lengths $|01|\cdot |02|\cdot |03|\cdot |04|$,
formula (\ref{F:M}) is equivalent to
\[
\sin \alpha\sin \gamma + \sin\beta\sin (\alpha + \beta + \gamma )
=
\sin (\alpha +\beta )\sin (\beta + \gamma ).
\]
But this is an identity, which follows easily by using the addition formula
for the sine function.
\qed

Note that the Gauss formula (\ref{F:G}) and the Monge formula (\ref{F:M}) are, in fact, equivalent.
Namely, if we write (\ref{F:G}) in the form $G=0$, and (\ref{F:M}) in the form $M=0$,
then it is easy to check that $G\equiv M$.

\begin{corollary}
\label{C:G}
Suppose that all vertex triangles of a convex pentagon have the
same area equal to some constant $k$. Then the area $A$ of the whole pentagon is
given by $A=\sqrt{5}\varphi k\approx 3.62 k$, where
$\varphi=\frac{1+\sqrt{5}}{2}$ is the golden ratio.
\end{corollary}

Corollary \ref{C:G} follows from (\ref{F:G}) by solving the equation $A^2-5kA+5k^2=0$.
In particular, for the regular pentagon of side length $a$ we get the well known value
$A=\frac{a^2}{4}\sqrt{25+10\sqrt{5}}\approx 1.72 a^2$.

As we shall prove in Section \ref{Sec:5}, a convex pentagon is affine--regular if and only if it has all vertex
triangles of the same area. Hence Corollary \ref{C:G} can be considered as a formula for the area of any
affine--regular pentagon.


Our next consequence generalizes the famous Ptolemy formula. In fact, the Monge
formula (\ref{F:M}) resembles the Ptolemy relation and since the sine--addition formula is
equivalent to the Ptolemy formula, no wonder that (\ref{F:M}) is a simple consequence of
the Ptolemy theorem. We leave this as a little challenge for the reader. But our next Corollary
implies the Ptolemy formula, and hence, finally, we can write symbolically:
\begin{center}
the Gauss pentagon formula $\Leftrightarrow$
the Monge formula
$\Leftrightarrow$
the Ptolemy formula.
\end{center}
\begin{corollary}
Let $0,1,2,3,4$ be any five points in a plane. Denote by $|ij|$ the distance between points $i$ and $j$
and by $R_{ijk}$ the circumradius of the triangle $\triangle ijk$, for $0\leq i,j,k\leq 4$. Then
\[
\frac{|12|}{R_{012}}
\cdot
\frac{|34|}{R_{034}}
+
\frac{|14|}{R_{014}}
\cdot
\frac{|23|}{R_{023}}
=
\frac{|13|}{R_{013}}
\cdot
\frac{|24|}{R_{024}}
\]
\end{corollary}
\proof
Just apply to every member of (\ref{F:M}) the formula for the area $A$
of a triangle : $A=\frac{abc}{4R}$. After dividing by $|01|\cdot|02|\cdot|03|\cdot|04|$, what is left is
the above equality. Later we shall prove that (\ref{F:M}) holds for any five points in a plane.
\qed

In the special case when all $R_{0ij}=R$, we get the Ptolemy formula.

\begin{corollary}
Let $C_1$ and $C_2$ be the first and the second cyclic symmetric functions
of the areas of border quadrilaterals of any pentagon with area $A$. Then
$A^2-C_1A+C_2=0$.
\end{corollary}
\proof Write the area of any quadrilateral as the difference between the
whole area $A$ and the area of the corresponding vertex triangle, e.g.\
$(1234)=A-(0)$, etc. This yields $C_1=5A-c_1$, and $C_2=5A^2-2c_1A+c_2$.
The claim now follows from Theorem \ref{T:Gpf}.\qed

\section{MORE GENERAL GAUSS PENTAGON FORMULA}
\label{sec:2}

First, let $T_i=(x_i)$, $i=0,1,2,3$ be any four points on a (coordinate) line $l$, and
$T_iT_j=x_j-x_i$ the oriented distance between $T_i$ and $T_j$. Then
\begin{equation}\label{E:1}
    T_1T_2\cdot T_3T_4+T_2T_3\cdot T_1T_4+T_3T_1\cdot T_2T_4 = 0.
\end{equation}
This follows from the following algebraic identity
\[
(x_2-x_1)(x_4-x_3)+(x_3-x_2)(x_4-x_1)+(x_1-x_3)(x_4-x_2)=0.
\]
Now take any point $O$ outside the line $l$ and let $h$ be its distance from the line $l$.
Denote by $(XYZ)$ the oriented area of the triangle $\triangle XYZ$. Then by multiplying
(\ref{E:1}) with $\frac{h^2}{4}$ and noting that $\frac{h}{2}T_iT_j=(OT_iT_j)=:t_{ij}'$,
it follows that (\ref{E:1}) is equivalent to
\begin{equation}\label{E:2}
    t_{12}'t_{34}'+t_{23}'t_{14}'+t_{31}'t_{24}'=0.
\end{equation}
Now let $P_1$, $P_2$, $P_3$, $P_4$ and $O$ be any five points in a plane.
Choose a line $l$ that intersect all four lines $OP_1$, $OP_2$, $OP_3$, $OP_4$ in
four points $T_1$, $T_2$, $T_3$, $T_4$, respectively.
\begin{figure}
\centerline{\includegraphics[width=131pt,height=141pt]{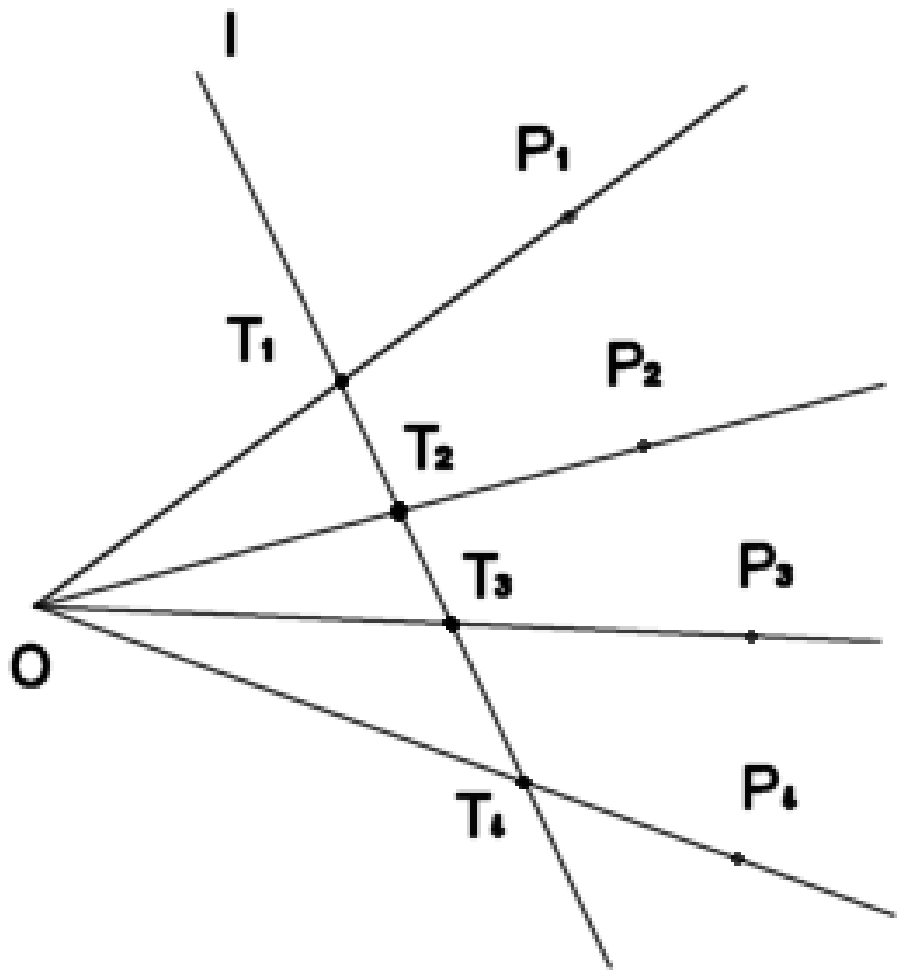}}
\caption{\label{Fig:2}}
\end{figure}
Let $t_{ij}=(OP_iP_j)$ and $t_{ij}'=(OT_iT_j)$ be oriented areas, and let $d$
be the distance from $T_i$ to $OP_j$ and $d'$ from $P_i$ to $OT_j$.
Then we have (for any $i\neq j$):
\[
\frac{t_{ij}}{t_{ij}'}=\frac{d\cdot |OP_j|}{d'\cdot |OT_j|}
=
\frac{|OP_i|}{|OT_i|}\frac{|OP_j|}{|OT_j|}=:k_i k_j.
\]
By using (\ref{E:2}) we then obtain
\[
t_{12}t_{34}+t_{23}t_{14}+t_{31}t_{24}
=
k_1k_2k_3k_4(t_{12}'t_{34}'+t_{23}'t_{14}'+t_{31}'t_{24}')=0,
\]
or writing again $(Oij)$ for $t_{ij}$ we get again the Monge formula (\ref{F:M}),
and even more general, namely for oriented triangle areas for any five
points in a plane.

Let us now reprove the Gauss formula. Again, let $P_1,\ldots, P_5$ be any five points
in a plane and $O$ any point in that plane, $t_{ij}=(Oij)$ the oriented area
of $\triangle OP_iP_j$, $t_{i}=(i)=(P_{i-1}P_iP_{i+1})$ area of the vertex triangle,
$A$ area of $P_1P_2P_3P_4P_5$, $A'$ area of $P_1P_3P_5P_2P_4$ (i.e.\ of the star pentagon),
$t_{i}'$ the oriented area of the vertex triangle of the star pentagon, e.g.\
$t_{1}'=(413)$, $t_{5}'=(352)$. Then clearly
\begin{equation}\label{E:3}
\begin{array}{r@{\ =\ }l}
A & t_{12}+t_{23}+t_{34}+t_{45}+t_{51},\\
A' & t_{13}+t_{35}+t_{52}+t_{24}+t_{41},\\
t_1 & t_{51}+t_{12}+t_{25} =t_{51}+t_{12}-t_{52},
\end{array}
\end{equation}
and four more similar equalities. From (\ref{E:3}) it follows easily that the first
and the second cyclic symmetric functions are given by
\begin{equation}\label{E:4}
\begin{array}{r@{\ =\ }l}
    c_1&\sum t_i=2A-A',\\
    c_2&\sum t_it_{i+1}=A(A-A').
\end{array}
\end{equation}
By permuting
$1\rightarrow 1$, $2\rightarrow 3$, $3\rightarrow 5$, $4\rightarrow 2$, $5\rightarrow 4$,
we see that $t_i$ becomes $t_{i}'$, $A$ becomes $A'$ and $A'$ becomes $-A$.
Hence the corresponding cyclic symmetric sums are as follows:
\begin{equation}\label{E:5}
\begin{array}{r@{\ =\ }l}
    c_1'&\sum t_i'=2A'+A,\\
    c_2'&\sum t_i't_{i+1}'=A'(A+A').
\end{array}
\end{equation}
From (\ref{E:4}) and (\ref{E:5}) we obtain the following generalized Gauss
pentagon type formulae.
\begin{theorem}\label{T:3.1}
\[
\begin{array}{l@{\ =\ }c@{\ =\ }l}
A^2-c_1A+c_2 & A^2-(2A-A')A+A(A-A') & 0,\\[2mm]
{A'}^2-c_1^2+4c_2 & {A'}^2-(2A-A')^2+4A(A-A') & 0,\\[2mm]
{A}^2-{c_1'}^2+4c_2' & {A}^2-(2A'+A)^2+4A'(A+A') & 0,\\[2mm]
{A'}^2-{c_1'}A'+c_2' & {A'}^2-(2A'+A)A'+A'(A+A') & 0.
\end{array}
\]
\end{theorem}

We can now give a geometric interpretation of the second solution of the Gauss formula (\ref{F:G}).
Namely, from the Equation (\ref{E:4}) it is evident that this solution is equal to $A-A'$, i.e.\ to the
difference between the areas of the original pentagon and its associated star pentagon. Similarly, from
the Equation (\ref{E:5}) we see that the second solution of the last equation in the Theorem \ref{T:3.1} is
$A+A'$.

\section{ANALOGUES OF GAUSS FORMULA FOR HEXAGONS, ETC.}
\label{sec:3}

Let $012345$ be a convex hexagon with area $A$, and write again $(ijk)$
for the area of the triangle $\triangle ijk$. For the vertex triangle
$\triangle (i-1,i,i+1)(\mbox{mod }6)$ we again simply write $(i)$ for its area.
Consider five pentagons having $0$ as a vertex and write Monge relation for
each of these pentagons in which all triangles have $0$ as a vertex.
Symbolically:
\begin{equation}
\label{S:1}
\begin{array}{lc@{\ =\ }c}
  \fbox{01234}: & (012)(034)+(014)(023) & (013)(024), \\[1mm]
  \fbox{01235}: & (012)(035)+(015)(023) & (013)(025), \\[1mm]
  \fbox{01245}: & (012)(045)+(015)(024) & (014)(025), \\[1mm]
  \fbox{01345}: & (013)(045)+(015)(034) & (014)(035), \\[1mm]
  \fbox{02345}: & (023)(045)+(025)(034) & (024)(035).
\end{array}
\end{equation}
\begin{figure}
\centerline{\includegraphics[width=142pt,height=122pt]{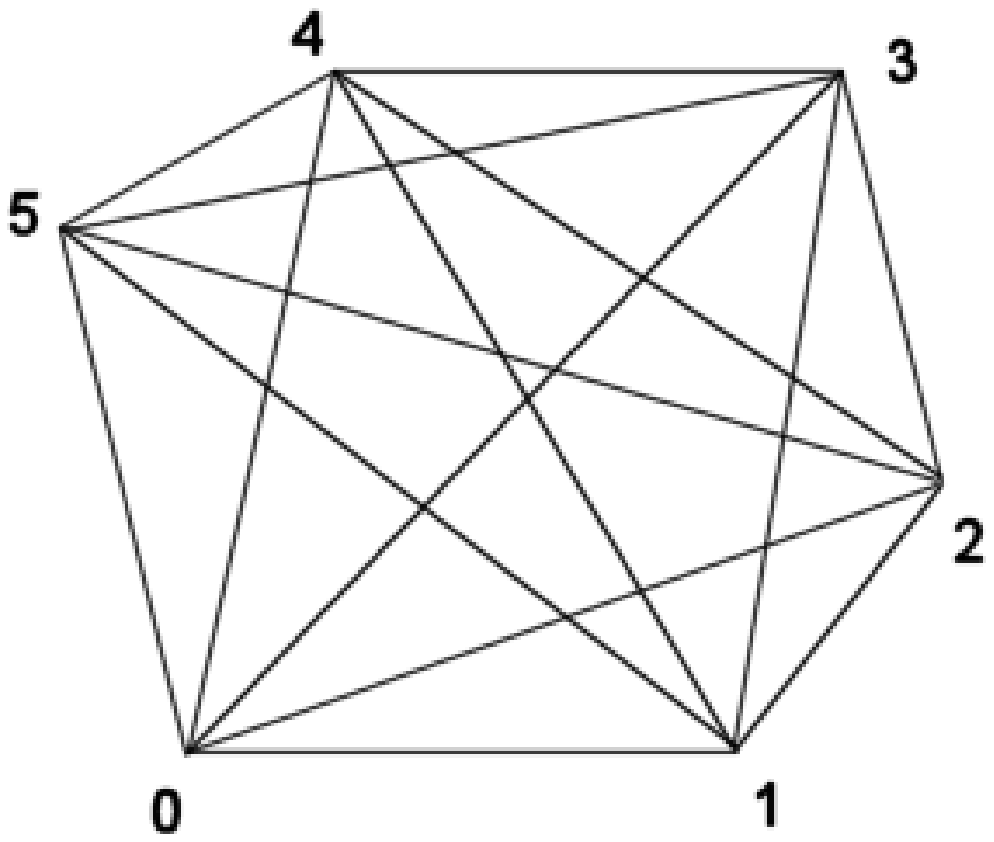}}
\caption{\label{Fig:3}}
\end{figure}
Certainly, here there is no formula for $A$ in terms of vertex triangles only. Therefore we choose an
additional parameter
\[
p=(013).
\]
By writing also $q=(014)$, $r=(023)$, $s=(025)$, it follows by using Figure \ref{Fig:3}
that we can make the following substitutions:
\[
\begin{array}{l}
(034)=A-({2})-({5})-p,\\
(023)=r=({2})+p-({1}),\\
(024)=A-({1})-({3})-({5}),\\
(035)=A-({2})-({4})-p,\\
(034)=A-({2})-({5})-p.
\end{array}
\]
Then the system (\ref{S:1}) becomes the following system of equations
in unknowns $A$, $q$ and $s$ with parameter $p$ :
\begin{equation}
\label{S:2}
\begin{array}{l}
  ({1})[A-(2)-(5)-p]+q[p+(2)-(1)]=p[A-(1)-(3)-(5)], \\[2mm]
  (1)[A-(2)-(4)-p]+(0)[p+(2)-(1)]=ps, \\[2mm]
  (1)(5)+(0)[A-({1})-({3})-({5})]=qs, \\[2mm]
  ({5})p+({0})[A-({2})-({5})-p]=q[A-({2})-({4})-p], \\[2mm]
  ({5})[p+({2})-({1})]+s[A-({2})-({5})-p]=\\
  \hspace{1cm}[A-({1})-({3})-({5})][A-({2})-({4})-p].
\end{array}
\end{equation}

It is not hard to check that the first three equations in (\ref{S:2}) imply the last two.

The solution of the first three equations gives again (as in the case of a pentagon)
a quadratic equation for $A$ in terms of all vertex triangles of the hexagon and one additional parameter
$p=(013)$. This solution obtained by {\tt Maple} gives the following hexagon analogue of the Gauss formula.
\begin{theorem}\label{T:4.1}
The area $A$ of a hexagon $012345$ with vertex triangle areas $(0)$, $(1)$, $(2)$, $(3)$, $(4)$, $(5)$ and
$p= $area$(013)$ satisfies the following equation.
\[
\begin{array}{l}
[p-(1)]A^2+[(1)(4)+2(1)(2)+(1)(5)-p^2+(1)p+(0)(1)-(2)p-\\
-(3)p-(5)p-(4)p-(0)(2)-(0)p]A-\\
-(1)(2)(5)-(1)(2)p-(1)(2)(4)-(1)(4)(5)-(1)(2)^2+(2)(3)p+(3)p^2-\\
-(0)(1)(5)+(0)(2)^2+(4)(5)p-(0)(1)(2)+(0)(2)(5)+(0)p^2-(0)(1)p+\\
+2(0)(2)p+(0)(5)p+(3)(4)p=0.
\end{array}
\]
\end{theorem}

It is easy to check that for $5\equiv 0$, this formula becomes the Gauss
formula~(\ref{F:G}).

\begin{corollary}
\label{C:1}
Suppose that all vertex triangles of a convex hexagon $012345$ have the same area equal to a constant $k$, and
let $p=$area$(013)=\lambda k$. Then the equation in $A$ from the Theorem \ref{T:4.1} reduces to the following equation
\[
(\lambda -1)A^2-(\lambda^2+4\lambda-4)kA+2(\lambda^2+2\lambda-2)k^2=0,
\]
or equivalently to
\begin{equation}\label{E:4.4}
A=\frac{\lambda^2+2\lambda-2}{\lambda -1}\ k
.\end{equation}
\end{corollary}
As we shall prove in the next section, a convex hexagon is affine--regular if and only if it has all vertex
triangles of the same area $k$ and one triangle, say $013$ (see the Figure \ref{Fig:3}) has the double area $2k$.
Hence, a hexagon is affine--regular if and only if all vertex triangles are equal to some constant $k$ and
$\lambda=2$. Then from (\ref{E:4.4}) we obtain $A=6k$, precisely what one would expect to obtain.

The following observation gives a rather symmetrical expression for $A$ which resembles the Gauss pentagon formula
(\ref{F:G}). Write (see Figure \ref{Fig:3})
\[
A=(135)+(0)+(2)+(4),\ \ \ A=(024)+(1)+(3)+(5),
\]
and multiply these two equations. Then we obtain
\begin{equation}\label{E:4.5}
\begin{array}{l}
A^2-[(0)+(1)+(2)+(3)+(4)+(5)]A+\\
+[(1)+(3)+(5)][(0)+(2)+(4)]=(135)(024).
\end{array}
\end{equation}
If we identify in the above equation (\ref{E:4.5}) the vertices $5$ and $0$, then the hexagon becomes a pentagon, and by using
the Monge formula (\ref{F:M}), we again infer the original Gauss pentagon formula (\ref{F:G}).

\begin{figure}
\centerline{\includegraphics[width=212pt,height=155pt]{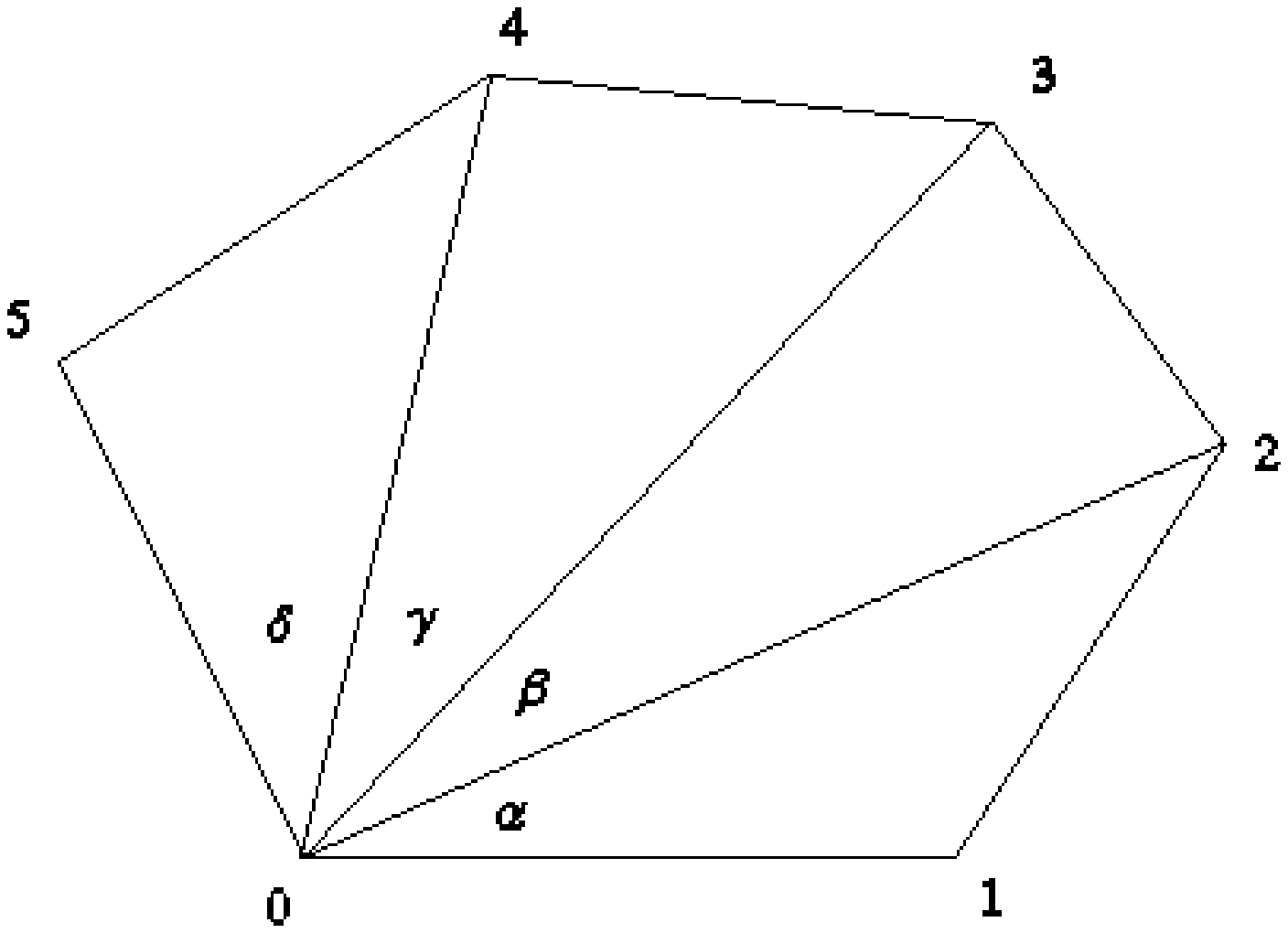}}
\caption{\label{Fig:4}}
\end{figure}

Another observation is the hexagon analogue of the Monge formula (\ref{F:M}):
\begin{equation}
\label{E:++}
(013)(024)(035)+(015)(023)(034)=(013)(025)(034)+(014)(023)(035).
\end{equation}
If one expresses the areas in terms of angles as on Figure \ref{Fig:4} and
in terms of side lengths $|0i|$, $1\leq i\leq 5$, then after dividing it by
$|01|\cdot |02|\cdot |03|\cdot |04|\cdot |05|$, what remains is an identity:
\[
\begin{array}{l}
\sin(\alpha+\beta)\sin(\beta+\gamma)\sin(\gamma+\delta)+
\sin(\alpha+\beta+\gamma+\delta)\sin\beta\sin\gamma=\\
=\sin(\alpha+\beta)\sin(\beta+\gamma+\delta)\sin\gamma+
\sin(\alpha+\beta+\gamma)\sin\beta\sin(\gamma+\delta).
\end{array}
\]
The formula (\ref{E:++}) (sometimes called Prouhet's formula) can easily
be obtained if we write it in the form $(013)[\cdots]=(023)[\cdots]$. Then
just using the last two equalities $\fbox{01345}$ and $\fbox{02345}$ from
the beginning of this section shows that it reduces to an obvious identity
\[
(013)(023)(045)=(023)(013)(045).
\]

Let us note the following fact for hexagons. Let $c_1$, $c_2$ resp.\
$C_1$, $C_2$ be the cyclic symmetric functions of areas of vertex
triangles resp.\ areas of border pentagons of any hexagon with area $A$.
Then $A^2-C_1A+C_2=A^2-c_1A+c_2$, which simply follows from $C_1=6A-c_1$,
$C_2=6A^2-2c_1A+c_2$.

In the case of a heptagon $0123456$, we would be tempted to apply the same
recipe as we did for hexagons. First to write down all $15$ Monge
relations for pentagons having $0$ as a vertex, then introducing two
parameters, say $p=(013)$ and $q=(014)$ and try to solve the obtained
system for the area $A$ of the heptagon in terms of areas of its vertex
triangles and parameters $p$ and $q$. We shall not consider it here.

Another way to naturally generalize the Gauss pentagon formula (\ref{F:G})
to any $n$--gon is to consider not only vertex triangles, but also
quadrilaterals of four consecutive vertices, five consecutive vertices and
so on. This cyclically recursive procedure could perhaps also be very
useful.

\section{AFFINE REGULAR PENTAGONS AND\\ HEXAGONS}
\label{Sec:5}

\begin{figure}
\centerline{\includegraphics[width=352pt,height=185pt]{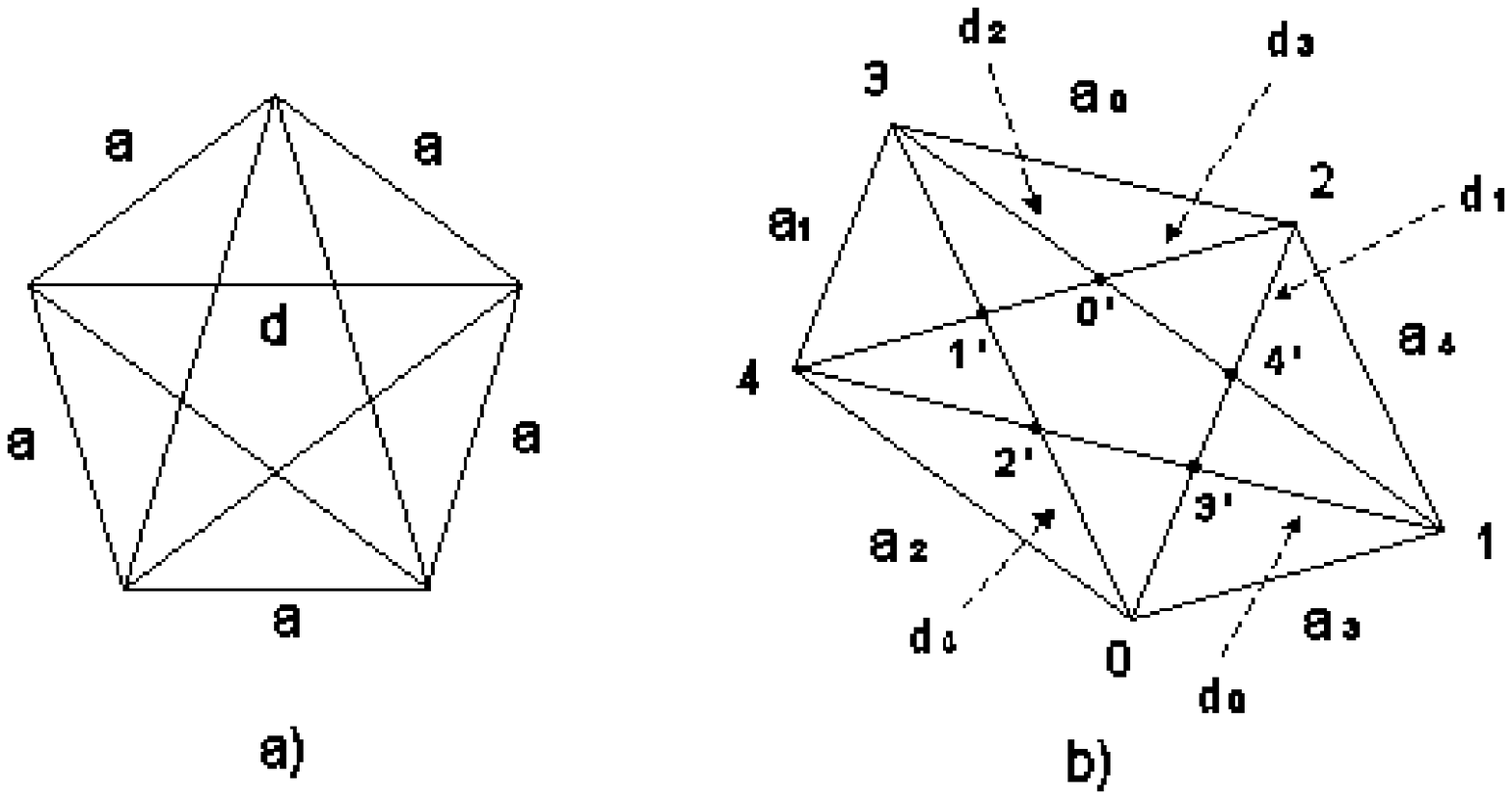}}
\caption{\label{Fig:5}}
\end{figure}

Recall that an {\it affine map} of the plane is a bijection of the plane
to itself such that any three collinear points are mapped again to three
collinear points. It maps lines to lines and parallel lines to parallel
lines. It preserves the ratio of (signed) distances between points on a
line. It also preserves the ratio of areas of two measurable sets in the
plane. An {\it affine regular polygon} is the image by an affine map of a
regular polygon. Clearly, any triangle is affine regular, and
parallelograms are precisely affine regular quadrilaterals. Any affine
regular $2n$--gon must be centrally symmetric, and hence its opposite
sides are parallel and have the same lengths. It is also obvious that any
affine regular polygon has all vertex triangles of the same area. We prove
that the converse holds for pentagons.
\begin{proposition}
\label{P:1}
A (convex) pentagon is affine regular if and only if it has all vertex triangles of the same area.
\end{proposition}
\proof
Let $01234$ be a convex pentagon that has all vertex triangles of the same area. The main feature of the
regular pentagon (Figure \ref{Fig:5}a) is that the ratio between its diagonal $d$ and the side $a$ is given
by
\[
\frac{d}{a}=\frac{a}{d-a}=\varphi,
\]
where $\varphi=(1+\sqrt{5})/2$ is the golden ratio.
By introducing notations as on Figure \ref{Fig:5}b, we argue as follows.
From the equality of areas of vertex triangles $(0)=(1)$, with the common base $a_3=|01|$,
it follows that vertices $2$ and $4$ are equally distant from $01$. So, $01||24$, i.e.\
$a_3||d_3$. Similarly $(1)=(2)$ implies $a_4||d_4$, and $(0)=(4)$ implies $a_2||d_2$.
Hence, we have two parallelograms $010'4$ and $0121'$. So, $|1'2|=|40'|=a_3$,
and $|1'0'|=2a_3-d_3$. A similar reasoning gives $|31'|=d_4-a_4$ and
$|30'|=d_2-a_2$. From the similar triangles $\triangle 1'0'3\sim\triangle 013$, we conclude
\[
\frac{d_4}{a_{4}}=\frac{d_2}{a_{2}}.
\]
In the same way we get $d_2:a_2=d_0:a_0$, etc. Finally, we obtain that the ratios of every diagonal $d_i$
and the corresponding side $a_{i}$ are equal. Call this ratio $\lambda$.
So, we have
\[
\frac{d_0}{a_0}=
\frac{d_1}{a_1}=
\frac{d_2}{a_2}=
\frac{d_3}{a_3}=
\frac{d_4}{a_4}=
\lambda\ .
\]
Once again, the similarity $\triangle 1'0'3\sim\triangle 013$ implies
\[
\frac{d_4-a_4}{d_4}=\frac{2a_3-d_3}{a_3},
\]
and so
\[
1-\frac{1}{\lambda}=2-\lambda \Rightarrow \lambda -\frac{1}{\lambda}=1 \Rightarrow
\lambda^2-\lambda-1=0 \Rightarrow \lambda=\varphi=\frac{1+\sqrt{5}}{2}.
\]
This is enough to conclude that $01234$ is affine regular, since all relevant ratios
are equal to the golden ratio, as it is in the regular pentagon.\qed

\begin{figure}
\centerline{\includegraphics[width=118pt,height=103pt]{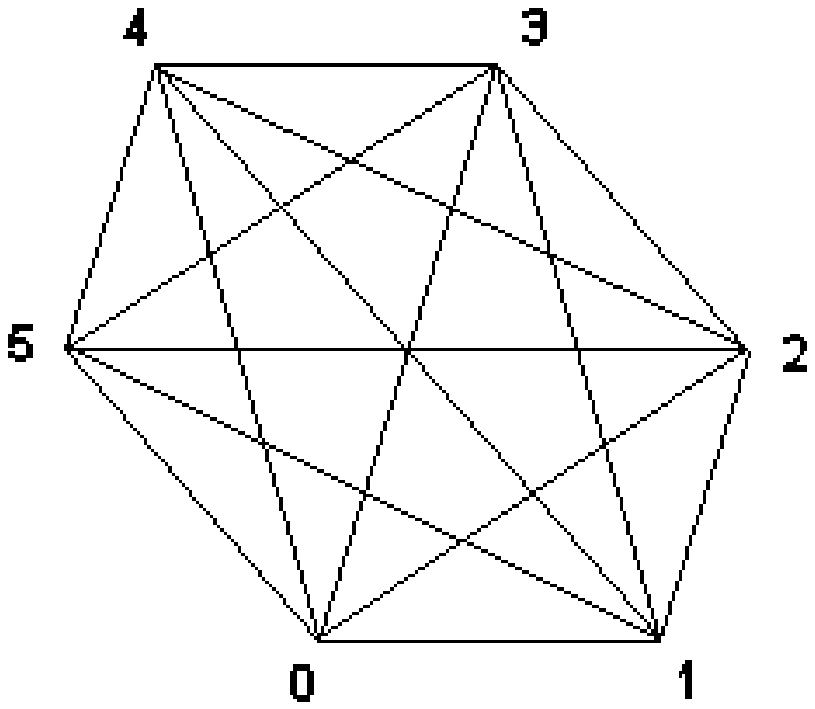}}
\caption{\label{Fig:6}}
\end{figure}

Now we give an analytic proof of a slightly more general claim.
\begin{proposition}
\label{P:2}
Let $A, B, C, D, E$ be any five different points in the plane, that are not all collinear.
Assume that
\begin{equation}
\label{F:||}\tag{$||$}
AB||EC,\ BC||AD,\ CD||BE,\ DE||CA,\ EA||DB.
\end{equation}
Then $ABCDE$ is an affine regular pentagon or a star affine regular pentagon.
\end{proposition}
\proof
Suppose some three of our five points are collinear. For example, let $A, B, C$ or $A, B, D$ are
collinear. From $BC||AD$ it would follow that $D$ or $C$ is on the same line.
So, all four points $A, B, C, D$ would be collinear. But, from $CD||BE$ we would conclude that $E$
also is on this line, what contradicts to our assumption. Hence, no three of our five points are collinear.
Therefore, there is a conic through these five points. By an affine transformation, if necessary,
we can assume that this conic is a parabola, or a rectangular hyperbola, or a circle. Hence, we have to examine
these three cases.
\begin{description}
  \item[$1^\circ$] If our conic is a parabola, we can assume that its equation is $y=x^2$ and let
  $a,b,c,d,e$ be abscissas of the points $A,B,C,D,E$ and we may assume $a<b<c<d<e$. The slope of the
  line $AB$ is equal to
  \[
\frac{a^2-b^2}{a-b}=a+b.
  \]
From the first two assumptions from (\ref{F:||}) we conclude that $a+b=e+c$, $b+c=a+d$. By adding,
we get $2b=e+d$ and this contradicts the facts that $b<e$ and $b<d$.

  \item[$2^\circ$] In the case of the rectangular hyperbola, we may take that its equation is $xy=1$.
  Under the same assumptions as in $1^\circ$, this time the slope of the line $AB$ is given by
  \[
\frac{\frac{1}{a}-\frac{1}{b}}{a-b}=-\frac{1}{ab}.
  \]
From (\ref{F:||}) it follows $ea=db$, $ab=ec$, $bc=ad$, $cd=be$, $de=ca$.
By multiplying two neighboring equalities we obtain (since $abcde\neq 0$):
$a^2=cd$, $b^2=de$, $c^2=ea$, $d^2=ab$. This implies that all five numbers have the same sign. Take, e.g.\
that they are all positive. If $a$ is the smallest, it follows $a<c$, $a<d$, and this contradicts $a^2=cd$.

  \item[$3^\circ$] It remains the case that our points $A,B,C,D,E$ are all on a circle. Then $BC||AD$
  implies that the lengths of chords $AB$ and $CD$ are equal, and similarly with the other assumptions of
  (\ref{F:||}). Hence, $ABCDE$ has all sides of the same length, and so it is either regular or star regular
  pentagon.\qed
\end{description}
\begin{proposition}
A (convex) hexagon is affine regular if and only if it has all vertex triangles of the same area and the area of
one of the triangles sharing only one side with the hexagon is twice the area of a vertex triangle.
\end{proposition}
\proof
Let again $012345$ be a hexagon with vertex triangles of equal areas, i.e.\ $(0)=(1)=\cdots=(5)$, and let
$(013)$ be twice $(0)$. See Figure \ref{Fig:6}.
From $(0)=(1)$ it follows that $25||01$, and from $(3)=(4)$ it follows also $25||34$. Hence $01||34$. Since
$(013)=2(1)$, it follows that the height of $\triangle 013$ is double than the height of the triangle $012$ on the
common base $\overline{01}$. Therefore, $25$ is equally distant from $01$ and $34$, and from $(1)=(3)$ it follows that
$|01|=|34|$. So, $0134$ is a parallelogram, and let $C$ be the intersection of its diagonals
$03$ and $14$. Now it is easy to see that the points of opposite sides of the hexagon $12$ and $45$ are equal and
parallel and also $05$ and $23$, and that $C$ is the center of symmetry of the hexagon. These facts easily imply
that our hexagon is affine regular.\qed

Note that results of this sections are in accordance with those in \cite{4}.

\section{THE ROBBINS FORMULA FOR A CYCLIC PENTAGON}

Let $01234$ be a (convex) {\it cyclic pentagon}, i.e.\ a pentagon which can be inscribed in a circle. Let $A$ be
the area of the pentagon, $R$ the circumradius, and $(i)=$ area$(\triangle i-1,i,i+1)$ (the indices taken modulo
$5$), $i=0,1,2,3,4$, the areas of vertex triangles. Denote by $a_{i}$ the side lengths of the pentagon and by
$d_i$ the lengths of the diagonals in the most systematic way as follows. Let $a_{i}$ be the length of the side
opposite to the vertex $i$, and let segments with lengths $a_{i}$ and $d_i$ be disjoint, $i=0,1,2,3,4$.
\begin{figure}
\centerline{\includegraphics[width=150pt,height=150pt]{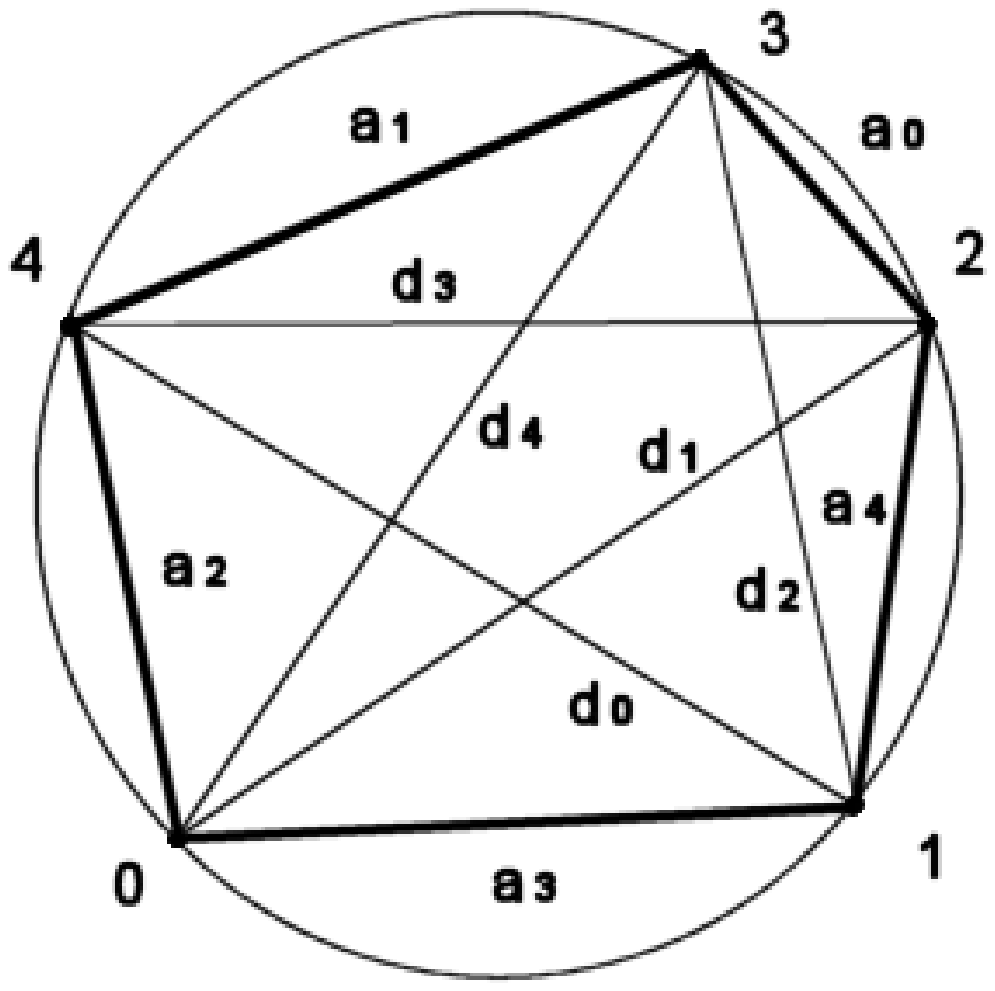}}
\caption{\label{Fig:7}}
\end{figure}
Our first aim is to find a formula for a diagonal in terms of the side lengths $a_{i}$'s. By cyclic symmetry it is
enough to find just one of them. Here is the result.
\begin{theorem}
\label{T:1} The diagonal $d_0$ (see Figure \ref{Fig:7}) satisfies the following polynomial equation of degree $7$ whose coefficients are
polynomials in $a_{i}$'s.
\begin{equation}\label{E:5.1}
\begin{array}{l}
\left(d_0^2-a_{2}^2-a_{3}^2\right)^2
(a_{0}d_0+a_{1}a_{4})
(a_{1}d_0+a_{0}a_{4})
(a_{4}d_0+a_{0}a_{1})
=\\
\hspace{2cm}=(a_{2}a_{3})^2
\left[
d_0^3-(a_{0}^2+a_{1}^2+a_{4}^2)d_0-2a_{0}a_{1}a_{4}
\right]^2.
\end{array}
\end{equation}
\end{theorem}
In other words, $X=d_0$ is one solution of the following degree $7$ equation:
\begin{equation}\label{E:5.2}
(X^2-q)^2(PX^3+SX^2+PQX+P^2)=p^2[X^3-QX-2P]^2,
\end{equation}
where
\begin{equation}\label{E:5.3}
\begin{array}{l}
p=a_{2}a_{3},\ P=a_{0}a_{1}a_{4},\ q=a_{2}^2+a_{3}^2,\ Q=a_{0}^2+a_{1}^2+a_{4}^2,\\[2mm]
S=(a_{0}a_{1})^2+(a_{0}a_{4})^2+(a_{1}a_{4})^2.
\end{array}
\end{equation}
Of course, there are four more similar formulas for diagonals $d_1, d_2, d_3, d_4$.
 In the special case $a_{0}=a_{1}=\cdots =a_{4}=a$, i.e.\ for the regular pentagon, it is very easy to check that
 $X=d_0=a\varphi$ is a solution, as expected ($\varphi$ is the golden ratio).

It can be proved that the associated polynomial of (\ref{E:5.1}) is irreducible over the rational functions
$\Q(a_{0},\ldots,a_{4})$. This is proved in \cite{10}.

To prove the theorem, we first need the following lemma about the (equivalent) claims for a
quadrilateral.

\begin{figure}
\centerline{\includegraphics[width=300pt,height=100pt]{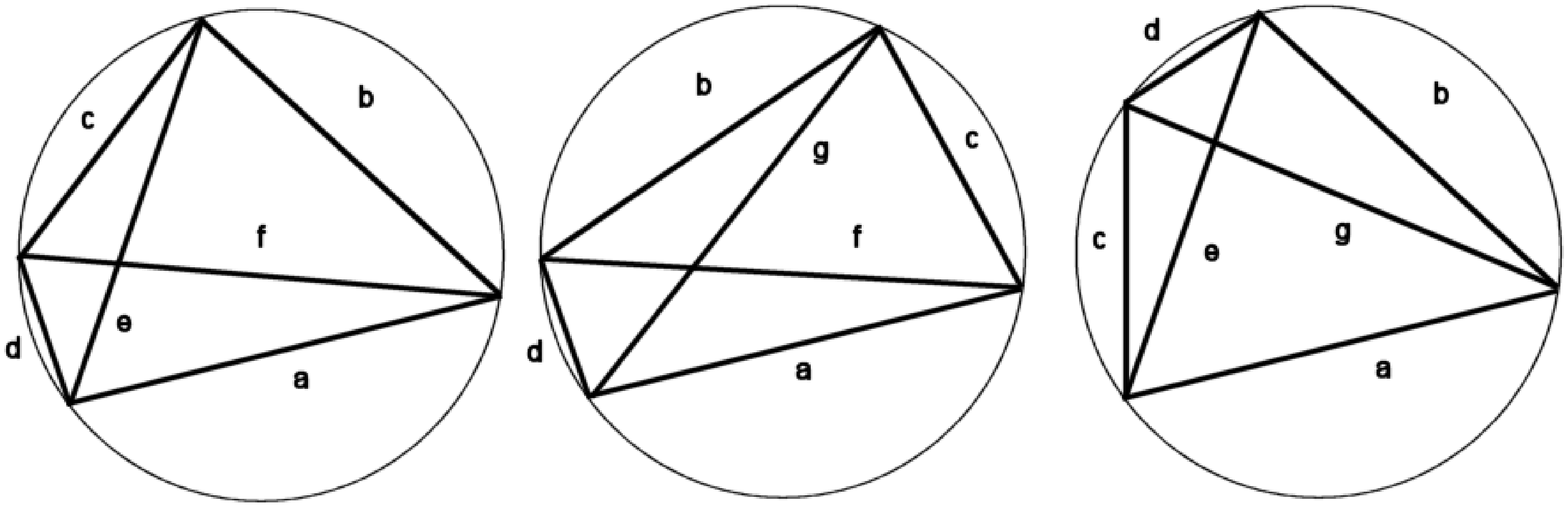}}
\caption{\label{Fig:8}}
\end{figure}

Consider three different cyclic quadrilaterals with side lengths $a, b, c, d$. They have the same area $A$,
the same perimeter $2s$ and the same circumradius $R$. Let $e, f, g$ be their diagonals (see Figure \ref{Fig:8}).

\begin{lemma}
\label{L:6.2}
For three cyclic quadrilaterals with the sides lengths $a, b, c, d$ diagonals $e, f, g$, area $A$,
perimeter $2s$ and circumradius $R$ the following relations hold.
\begin{align}
a)\  & 4AR=(ab+cd)e,\ 4AR=(ad+bc)f,\ 4AR=(ac+bd)g;\label{Eq:6.14} \\
b)\  & 4AR=efg;\notag  \\
c)\  & ef=ac+bd,\ fg=ab+cd,\ eg=ad+bc\ (\mbox{Ptolemy formula, 2nd century});\notag\\
d)\  & (4AR)^2=(ab+cd)(ac+bd)(ad+bc),\label{Eq:6.15}\\    &
\mbox{or expanded by $a$}:\notag\\    &
(4AR)^2=Pa^3+Sa^2+PQa+P^2,\label{Eq:6.16}\\    & \mbox{where } P=bcd,\
S=(bc)^2+(bd)^2+(cd)^2,\ Q=b^2+c^2+d^2;\notag\\
e)\  & e^2=\frac{(ad+bc)(ac+bd)}{ab+cd}, f^2=\frac{(ab+cd)(ac+bd)}{ad+bc},\notag\\    &
g^2=\frac{(ab+cd)(ad+bc)}{ac+bd};\label{Eq:6.17}\\
f)\  & (4R)^2=\frac{(ab+cd)(ac+bd)(ad+bc)}{(s-a)(s-b)(s-c)(s-d)},\label{Eq:6.18}\\
g)\  & (4A)^2=2[(ab)^2+(ac)^2+(ad)^2+(bc)^2+(bd)^2+(cd)^2]-\notag\\    &
-(a^4+b^4+c^4+d^4)+8abcd=16(s-a)(s-b)(s-c)(s-d)\label{Eq:6.19}\\    &
\mbox{(Brahmagupta's formula, $7$th century).}\notag
\end{align}
\end{lemma}
\proof
A short proof in the classic geometric combinatorial style is as follows.\\[2mm]
\noindent
$a)$ Apply the well--known triangle formula $4\Delta R=abc$ to the six triangles separated by diagonals $e, f, g$,
respectively, as on Figure \ref{Fig:8}.\\
$b)$ and $c)$ By equating the equalities from $a)$ and dividing by $efg$ we obtain
\[
\frac{ab+cd}{fg}=\frac{ad+bc}{eg}=\frac{ac+bd}{ef}=\frac{4AR}{efg}=:k,
\]
where $k$ is a constant. But, when $d=0$, then $e=c$, $f=a$, $g=b$ (see Figure \ref{Fig:8}) and so $k=1$.\\
$d)$ Just multiply the equalities from $c)$ and use $b)$.\\
$e)$ Multiply any two of the three equalities from $c)$ and divide by the third one.\\
$f)$ Write the formula for the area of the triangle with sides lengths $a,
b, e$ (see (\ref{E:6.21})) in the form
\[
R^2(2a^2b^2+2a^2e^2+2b^2e^2-a^4-b^4-e^4)=a^2b^2e^2,
\]
Then substitute $e^2$ from the formula $e)$ and simplify.\\
$g)$ Divide the results $d)$ and $f)$.
\qed

We used here the famous Heron formula (from the 1st century, which actually goes back to Archimedes)
for the area $\Delta$ of the triangle with side lengths $a, b, c$:
\begin{equation}\label{E:6.21}
    (4\Delta)^2=2[(ab)^2+(bc)^2+(ca)^2]-(a^4+b^4+c^4)=:H^2(a^2,b^2,c^2).
\end{equation}
Also, write $B^2(a^2,b^2,c^2,d^2)$ for the right hand side of
Brahmagupta's formula (\ref{Eq:6.19}). Note (see Figure \ref{Fig:8}):
\[
B^2(a^2,b^2,c^2,d^2)=[H(a^2,b^2,e^2)+H(c^2,d^2,e^2)][H(a^2,d^2,f^2)+H(b^2,c^2,f^2)],
\]
where $e, f$ are given by
(\ref{Eq:6.17}).

\proof[of Theorem \ref{T:1}] Dissect the pentagon $01234$ with diagonal of length $d_0$ into the vertex triangle
$401$ with area $(0)$ and the remaining quadrilateral $1234$. Then
\begin{equation}\label{E:6.22}
    A=(0)+area(1234).
\end{equation}
Then for the radius $R$ applied to the vertex triangle and to the quadrilateral, using $d)$ from
the Lemma \ref{L:6.2}, we have
\begin{equation}\label{E:6.23}
    4(0)R=a_{2}a_{3}d_0,
\end{equation}
\begin{equation}\label{E:6.24}
    [4area(1234)R]^2=(a_{0}d_0+a_{1}a_{4})(a_{1}d_0+a_{0}a_{4})(a_{4}d_0+a_{0}a_{1}).
\end{equation}
Since $[4(0)]^2=H^2(a_{2},a_{3},d_0)$, and $[4area(1234)R]^2=B^2(a_{0},a_{1},a_{4},d_0)$,
if we write $X=d_0$, then with notations from (\ref{E:5.3})we have
\begin{equation}\label{E:6.25}
    H^2(a_{2},a_{3},X)=4p^2-(X^2-q)^2,
\end{equation}
\begin{equation}\label{E:6.26}
B^2(a_{0},a_{1},a_{4},X)=4(S+2PX)-(X^2-Q)^2,
\end{equation}
and
\begin{equation}\label{E:6.27}
(a_{0}X+a_{1}a_{4})(a_{1}X+a_{0}a_{4})(a_{4}X+a_{0}a_{1})=PX^3+SX^2+PQX+P^2.
\end{equation}
By (\ref{E:6.23}) and (\ref{E:6.24}) we have
\begin{equation}\label{E:6.28}
    [4p^2-(X^2-q)^2]R^2=p^2X^2
\end{equation}
and
\begin{equation}\label{E:6.29}
    [4(S+2PX)-(X^2-Q)^2]R^2=PX^3+SX^2+PQX+P^2.
\end{equation}
By simply eliminating $R^2$ from (\ref{E:6.28}) and (\ref{E:6.29}) we get (\ref{E:5.2}).\qed

Note that if we put $a_{4}=0$ in (\ref{E:5.1}), then it is easy to see that it collapses
to the formula for a diagonal of a cyclic quadrilaterals in terms of the side lengths, i.e.\ to
(\ref{Eq:6.17}).

\begin{proposition}
\label{P:6.3}
The diagonal $d_0=X$ of the cyclic pentagon in terms of the side lengths $a_{i}$'s as on the
Figure \ref{Fig:7}, the area $A$, and the circumradius $R$ is given by the following cubic
(in notations of (\ref{E:5.3})):
\begin{equation}\label{E:6.30}
    (4AR-pX)^2=PX^3+SX^2+PQX+P^2
\end{equation}
\end{proposition}
\proof
This follows again by dissecting the pentagon into the vertex triangle $401$ and the remaining
quadrilateral and using Lemma \ref{L:6.2} $d)$.\qed
\begin{corollary}
With notations as in Proposition \ref{P:6.3} we have the following quadratic:
\begin{equation}\label{E:6.31}
    (X^2-q)4AR=p[(Q-q)X+2P].
\end{equation}
\end{corollary}
\proof If we put the right hand side of (\ref{E:6.30}) into (\ref{E:5.2}) we get:
\[
(X^2-q)^2(4AR-pX)^2=p^2(X^3-QX-2P)^2.
\]
Hence we obtain
\[
(X^2-q)(4AR-pX)=\pm p(X^3-QX-2P).
\]
The sign $+$ leads to a contradiction, while with the sign $-$ we have the quadratic polynomial
(\ref{E:6.31}).\qed

Our next aim is to eliminate $X$ from (\ref{E:6.30}), (\ref{E:6.31}) and (\ref{E:6.28})
to get various relations between $R, A$ and $p, q, P, Q, S$.
These parameters are related to symmetric functions of $a_{i}$'s as follows.

Let $a_{i}^2=x_i$, $i=0,1,2,3,4$ and let $e_k=e_k(x_0,x_1,x_2,x_3,x_4)$, $k=1,2,3,4,5$
be the elementary symmetric functions (see \cite{9}). Then it is easy to check:
\begin{equation}\label{E:6.33}
\begin{array}{c}
q+Q={e}_1, S+p^2+qQ={e}_2, qS+p^2Q+P^2={e}_3,\\
p^2S+P^2q={e}_4, (pP)^2={e}_5.
\end{array}
\end{equation}

By writing $H=H(a_{2}^2,a_{3}^2,X^2)$, $B=B(a_{0}^2,a_{1}^2,a_{4}^2,X^2)$, where $X=d_0$, for the area $A$
of the pentagon we have $(4A)^2=(H+B)^2=H^2+B^2+2HB$, or $[(4A)^2-(H^2+B^2)]^2=4H^2B^2$. By substituting
(\ref{E:6.25}) and (\ref{E:6.26}) we get the following degree $7$ polynomial equation in X:
\begin{equation}\label{E:6.33a}
\begin{array}{l}
\{(4A)^2-[4p^2-(X^2-q)^2+4(S+2PX)-(X^2-Q)^2]\}=\\
\hspace{6.2ex}=4[4p^2-(X^2-q)^2][4(S+2PX)-(X^2-Q)^2].
\end{array}
\end{equation}

First, we relate $4AR$ and side lengths.

Let us write the formula (\ref{Eq:6.15}) in a more appropriate way in terms of the symmetric function theory.
Namely, write (\ref{Eq:6.15}) in terms of the elementary symmetric functions $e_k=e_k(a^2,b^2,c^2,d^2)$,
$k=1,2,3,4$.

See \cite{9} for definitions and properties. Write (\ref{Eq:6.15}) for $Z=(4AR)^2$ in the form
\begin{equation}\label{E:6.34}
\begin{array}{l}
    Z=(4AR)^2=\\
    =a^3bcd+ab^3cd+abc^3d+abcd^3+\\
    \hphantom{=}+a^2b^2c^2+a^2b^2d^2+a^2c^2d^2+b^2c^2d^2=\\
    =e_1\sqrt{e_4}+e_3.
\end{array}
\end{equation}


In the form of polynomials in $e_k$'s, this can be written as
\begin{equation}\label{Eq:6.34}
    (Z-e_3)^2-e_1^2e_4=0.
\end{equation}
The following is the pentagon analogue of (\ref{Eq:6.34}).
\begin{proposition}
\label{P:6.5}
Let $A$ be the area, $R$ the radius of a cyclic pentagon, and $e_k$, $1\leq k\leq 5$, the elementary symmetric
functions in the squares of the side lengths. Then $Z=(4AR)^2$ satisfies the following degree $7$ polynomial equation.
\begin{equation}\label{Eq:6.35}
\begin{array}{l}
  Z^3[(Z-e_3)^2-e_{41^2}]^2+e_5[-e_{51^3}+5e_{521}-8e_{53}+(-e_{52^21^4}-\\
  -4e_{5321^3}+8e_{52^31^2}-32e_{5421^2}-4e_{53^21^2}+16e_{532^21}-16e_{52^4}-\\
  -64e_{5431}+128e_{542^2}-256e_{54^2}-e_{31^2}-8e_{41}+4e_{32}-24e_5)Z+\\
  +(2e_{421^4}+28e_{521^3}+4e_{31^3}-8e_{42^21^2}-2e_{3^221^2}+56e_{531^2}-112e_{52^21}+\\
  +32e_{4^21^2}-4e_{3^31}+8e_{3^22^2}+448e_{541}-32e_{43^2}+e_{1^2}+4e_2)Z^2+\\
  (-28e_{41^3}+4e_{321^2}-196e_{51^2}+36e_{3^21}-16e_{32^2}+64e_{43})Z^3+\\
  +(-2e_{21^2}-60e_{31}+8e_{2^2}-32e_4)Z^4+28e_1Z^5]=0.
\end{array}
\end{equation}
\end{proposition}
Here we used the standard notation (see \cite{9}) $e_\lambda=e_{\lambda_1}e_{\lambda_2}\ldots$ for a partition
$\lambda=(\lambda_1,\lambda_2\ldots)=(n^{m_n}\ldots 3^{m_3}2^{m_2}1^{m_1})$,
where $m_i=\#\{j|\lambda_j=i\}$.

\proof
This follows by eliminating $X$ via resultant (see \cite{8}) from the cubic
(\ref{E:6.30}) and the quadratic (\ref{E:6.31}). Then we express the result in
terms of monomial symmetric functions of $a_i^2$'s and by using Stembridge's
{\it SF} {\tt Maple} package (\cite{01}) we get the result in terms of elementary
symmetric functions.
\qed

If we let one side length to be zero, e.g.\ $a_{4}=0$, then the pentagon collapses to a quadrilateral. It is clear
that in that case $e_5=0$, and (\ref{Eq:6.35}) simply reduces to (\ref{Eq:6.34}).

Similarly, the formula (\ref{Eq:6.18}) for the radius $R$ has the following pentagon analogue.
\begin{proposition}
\label{P:6.6}
The square of the circumradius $R$ of a cyclic pentagon in terms of the elementary symmetric functions $e_k$
in the squares of the sides satisfies the following degree $7$ equation (in $R^2$):
\begin{equation}\label{Eq:6.36}
\begin{array}{l}
    R^6[(e_{1^4}-8e_{21^2}+16e_{2^2}-64e_4)R^4+(2e_{31^2}+16e_{41}-8e_{32})R^2-\\
    -e_{41^2}+e_{3^2}]^2+e_5W=0,
\end{array}
\end{equation}
where $W$ equals
\[
\!\!\begin{array}{l}
2048(-e_{1^3}+4e_{21}-8e_3)R^{14}+32(23e_{1^4}-88e_{21^2}+192e_{31}-16e_{2^2}+64e_4)R^{12}\\
+64(-e_{1^5}+2e_{21^3}-9e_{31^2}+8e_{2^21}-8e_{41}-12e_{32}-12e_5)R^{10}\\
(e_{1^6}+6e_{21^4}+32e_{31^3}-32e_{2^21^2}-32e_{41^2}-32e_{2^3}+256e_{51}+128e_{42}+224e_{3^2})R^8\\
+2(-e_{31^4}+4e_{41^3}+2e_{321^2}-8e_{51^2}-16e_{3^21}+8e_{32^2}-16e_{52}-32e_{43})R^6\\
(2e_{51^3}-2e_{421^2}+e_{3^21^2}+16e_{531}-8e_{52^2}+8e_{431}-2e_{3^22}-16e_{54}+8e_{53})R^4\\
+(-2e_{31}+e_{2^2}-4e_4)R^2e_5+e_{5^2}.
\end{array}
\]
\end{proposition}
Finally, by eliminating $X$ and $R$ from the quadratic , cubic and quartic
equations (\ref{E:6.31}), (\ref{E:6.30}) and (\ref{E:6.28}), and by changing the
basis using (\ref{E:6.33}), we get by using {\tt Mathematica} or {\tt Maple} the
following main result. \begin{theorem}
\label{T:6.7}
{\bf (the Robbins formula)} Let $A$ be the area of a cyclic pentagon and $e_k$, $1\leq k\leq 5$, the elementary symmetric functions
in the squares of the side lengths. Then $Y=(4A)^2$ satisfies the following monic degree $7$ polynomial
equation
\begin{equation}\label{E:**}\tag{**}
\begin{array}{l}
[(Y-4e_2+e_1^2)^2-64e_4]^2\{Y[(Y-4e_2+e_1^2)^2-64e_4]+\\[1mm]
[e_1(Y-4e_2+e_1^2)+8e_3]^2\}-128e_5\{16[e_1(Y-4e_2+e_1^2)+8e_3]^3+\\[1mm]
18Y[e_1(Y-4e_2+e_1^2)+8e_3]+[(Y-4e_2+e_1^2)^2-64e_4]+2^73^3Y^2e_5\}=0.
\end{array}
\end{equation}
\end{theorem}
We remark that the coefficients of powers of $Y$ in (\ref{E:**}),
considered as polynomials in the elementary symmetric functions of the
sides length squared are primitive (i.e.\ their integer coefficients are
relatively prime). It is easy to check that our formula agrees with the
Robbins formula (\cite{6}, \cite{77}) and hence it can be considered as a
new proof of Robbins formula. Again, if we let one side length to be zero,
e.g.\ $a_{4}=0$, then the pentagon collapses to a quadrilateral and then
$e_5=0$. In that case (\ref{E:**}) reduces to the Brahmagupta formula
(\ref{Eq:6.19}) written for $Y=(4A)^2$ in term of $e_k$'s as follows:
\[
(Y-4e_2+e_1^2)^2-64e_4=0.
\]

By combining Gauss and Robbins formulas for cyclic pentagon we get the following
unexpected result.
\begin{theorem}
\label{T:6.8}
The area of a cyclic pentagon is a rational function in areas of vertex triangles
and side lengths squared.
\end{theorem}
\proof
By using the Gauss pentagon formula of Theorem \ref{T:Gpf}, written as
\[
(4A)^2=c_1'\cdot 4A+c_2',\ (\mbox{with }\ c_1'=4c_1\ \mbox{ and }\ c_2'=16c_2)
\]
as a side relation, we can simplify the Robbins formula of Theorem ref{T:6.7}
for $Y=(4A)^2$:
\[
  Y^7+C_1Y^6+C_2Y^5+C_3Y^4+C_4Y^3+C_5Y^2+C_6Y+C_7=0,
\]
to a linear equation for $A$ in  terms of the cyclic invariants $c_1', c_2'$
of the areas of vertex triangles and symmetric functions $C_1,\ldots, C_7$,
given by Theorem \ref{T:6.7}, of
side lengths squared. The result is the following formula for the area of a
cyclic pentagon
\[
A=\dfrac{N}{D},
\]
where (with $v={c_1'}^2-{c_2'}$ and $u={c_1'}^2-2{c_2'}$)
\[
\begin{array}{l}
N=
{c_2'}[(v^5-4v^4{c_2'}+2v^3{c_2'}^2+5v^2{c_2'}^3-2v{c_2'}^4-{c_2'}^5)C_1
+v(v^3-3v^2{c_2'}+3{c_2'}^3)C_2\\
+(v^3-2v^2{c_2'}-v{c_2'}^2+{c_2'}^3)C_3
+(v^2-v{c_2'}-{c_2'}^2)C_4+vC_5\\
+v^6-5v^5{c_2'}+5v^4{c_2'}^2+6v^3{c_2'}^3-7v^2{c_2'}^4-2v{c_2'}^5+{c_2'}^6+C_6]-C_7
\end{array}
\]
and
\[
\begin{array}{l}
D=
c'_1[u(u^2-{c_2'}^2)(u^2-3{c_2'}^2)C_1
+((u^2-{c_2'}^2)^2-u^2{c_2'}^2)C_2
+u(u^2-2{c_2'}^2)C_3\\
+(u^2-{c_2'}^2)C_4
+uC_5
+C_6
+u^2(u^2-2{c_2'}^2)^2-{c_2'}^2(u^2-{c_2'}^2)^2].\qed
\end{array}
\]

By a similar procedure we can also get the area of a cyclic hexagon and confirm one more result of
Robbins (\cite{6}, \cite{77}), but we shall not present it here. Let us mention only that the
computational complexity grows very rapidly, as Robbins himself remarked also about ten years ago.
More on that in \cite{7} and \cite{00}.

We also mention that (\ref{E:**}) as well as (\ref{Eq:6.35}) and (\ref{Eq:6.36})
can be written in other base for symmetric functions: monomial functions, power
sums, Schur functions (for basics on symmetric functions see \cite{9}), but
elementary symmetric functions provide the most compact and shortest way to
write down these formulas.


\section{SOME FINAL REMARKS AND PROBLEMS}

\ \\
\indent {\bf 1.} Of course, a recursive procedure of our method could, at least in principle, be applied to any cyclic $n$--gon.
In fact, it can be proved inductively that the length of any diagonal of a cyclic $n$--gon
is an algebraic function of the side lengths. The same applies to the circumradius $R$
and the area $A$. These facts are proved in \cite{10}.


{\bf 2.} By using rigidity theory and algebraic geometry, it was quite recently proved in \cite{11} (not published yet)
that in general, for any cyclic $n$--gon there is a (nonzero) polynomial equation in $(4A)^2$ whose coefficients
$C_t$ are polynomials in the side lengths squared. The smallest degree of such a polynomial is equal to $\Delta_k$
for $n=2k+1$ and $2\Delta_k$ for $n=2k+2$, where
\[
\Delta_k=\frac{1}{2}[(2k+1)
\left(\!
\begin{array}{c}
2k\\k
\end{array}\!
\right)-2^{2k}
]=\sum_{i=0}^{k-1}(k-i)
\left(\!
\begin{array}{c}
2k+1\\i
\end{array}\!
\right)
\]
is the total number of $(2k+1)$--gons (convex or nonconvex) inscribed in a
circle. This was also conjectured by Robbins in \cite{6}, \cite{77}. It is
also recently proved by R.~Connelly \cite{13} that these polynomials are
monic. We conjecture the following.

\begin{conjecture}
For any cyclic polygon with area $A$ the coefficients of the polynomial
equation for $(4A)^2$, considered as polynomials in the elementary
symmetric functions of the sides length squared, are primitive (i.e.\
their integer coefficients are relatively prime).
\end{conjecture}

{\bf 3.} It should be emphasized that the coefficients $C_t$ in the Robbins formula
are symmetric polynomials in side length squared.
But for general pentagon, however, the Gauss quadratic formula (\ref{F:G}) for $A$
is invariant only with respect to cyclic symmetries.
By combining the Gauss and Robbins formula, one can expect many new interesting
problems and facts (like Theorem \ref{T:6.8}) in geometry of pentagons.

Another one. How to construct by a ruler and compass a cyclic pentagon with given
side lengths? Well, we have the $7$--th degree irreducible equation for a diagonal
(or radius $R$). Hence, in general, it is not possible to do it. Some questions are
also posed in \cite{13}, \cite{11}, \cite{12}, \cite{14}, \cite{10}.

So, after centuries of geometry of mostly  triangles and quadrilaterals, we finally arrived to the (nontrivial)
geometry of pentagons.

More about $5$--gons, $6$--gons etc.\ and some computations will appear elsewhere
(\cite{7}, \cite{00}).

{\bf 4.} We also remark here that it seems possible that our method could
as well be applied for computing volumes of polytopes in any dimension.
The central role here would be the general Pl\" ucker formulas (see
\cite{8}) which reduce to Monge relation in dimension two. And even
obtained formula for areas contain some information about volumes, since,
for example, any convex pentagon is a projection of a $4$--simplex, or of
a $3$--dimensional bipyramid to the plane etc.

{\bf 5.} One can also hope to use these ideas for computing volumes of at
least simplicial polytopes or polytopes whose vertex figures are
simplices. Then one would use the Heron formula for simplices (i.e.\
Cayley--Menger's formula) etc. Ideas in \cite{12} and \cite{13} can
certainly be fruitful here.

To end this general nonsense, one could confirm once again that, no matter
how, but computing volume is hard, indeed.
\vspace{5mm}

\begin{acknowledgement}
We thank Professor Victor
\mbox{Adamchik} from Carnegie Mellon University whose expertise in
computer algebra helped us essentially in obtaining our results while he was
visiting Zagreb University. We also thank Dr. Igor Urbiha for his help in
writing this paper.
\end{acknowledgement}

\end{document}